\documentclass[12pt]{article}%
\usepackage{amsfonts}
\usepackage{amsmath,amscd,amsthm,a4,amssymb}
\usepackage{amsmath,amscd,amsthm,a4,amssymb}
\usepackage{amsmath}
\usepackage{amssymb}
\usepackage{graphicx,multirow,rotating}%

\providecommand{\U}[1]{\protect\rule{.1in}{.1in}}

{}

{}

{}

{}

{}
\theoremstyle{plain}
{}

{}

{}

{}

{}

{}

{}

{}

{}

{}

{}

{}

{}

{}

{}

{}
\begin{document}

\begin{center}
{\Large \textbf{A New Descent Algebra of Weyl Groups of Type A}}$_{n}%
${\Large \textbf{\ \ }}
\end{center}

\centerline{\large T\"{u}lay YA\v{G}MUR $^{1}$, Himmet CAN $^{2}${\footnotetext{
{E-mail: $^{1}$tyagmur@erciyes.edu.tr (T.Ya\v{g}mur); $^{2}$can@erciyes.edu.tr (H.Can)}} }}

\ 

\centerline{\it $^{1,2}$ Department of Mathematics, Erciyes University, {\textsc{38039}}
Kayseri, Turkey}

\begin{abstract}
In this paper we define an equivalence relation on the set of all
$x_{J}$ in order to form a basis for a new descent algebra of Weyl groups of type  $A_{n}$. By means of this,  we construct a new commutative and semi-simple  descent algebra of Weyl groups of type $A_{n}$   generated
by equivalence classes arising from this equivalence relation.

\end{abstract}

\begin{quote}

\end{quote}
 \bigskip
\section{Introduction}

\bigskip
 The main objective of this paper is to construct a new descent algebra of Weyl groups of type $A_{n}$ by using Solomon's basic concepts.  Moreover,  the basis of this new algebra consists of equivalence classes arising from the equivalence relation on the set $x_J$, which is going to be defined in this section. This new descent  algebra  will be denoted by $\sum_{W}(A_{n})$.\\
\indent  In the next section, we shall determine two basic properties of this new descent algebra. That is to say, by contrast with the Solomon's  descent algebra, we shall show that this new algebra is commutative and also semi-simple. 
Furthermore, we shall show that there is an  isomorphism between the new descent algebra $\sum_{W}(A_{n})$ and the parabolic Burnside ring of associated Weyl group. Note that, there is a homomorphism between the Solomon algebra $\sum{(W)}$  and the parabolic Burnside ring $PB(W)$ [3]. In addition, we shall show that there is an isomorphism between the new descent  algebra $\sum_{W}(A_{n})$ and the $\mathbb{Z}$-module of class functions $G(W(A_n))$ generated by permutation character ${\chi}_J$.\\
\indent Now the notation, which is fairly standard and follows that given in Carter [1], and  Solomon [2] or Bergeron et al. [3], is  first summarized.\\ 
\indent Let $\Phi$ be a root system with fundamental system $\Pi$ and corresponding positive system ${\Phi}^{+}$. If  $r\in{\Phi}$ , let $w_r$ denote the reflection corresponding to $r$. Besides, let $W=W({\Phi})=W({\Pi})$ be the {\textit{ Weyl group}} of root system $\Phi$, that is, the group generated by the $w_r$, $r\in{\Phi}$.\\ \\
\indent Let  $e_{1}, e_{2}, ... , e_{n+1}$   be an orthonormal basis of a Euclidean space of dimension $n+1$.
Then the fundamental system and root system of type $A_{n}$ are given by 
\begin{center}
 $\Pi =\{e_{1}-e_{2}, e_{2}-e_{3}, ...  , e_{n-1}-e_{n}, e_{n}-e_{n+1}\}$, 
\end{center}

\begin{center}
 $\Phi =\{e_{i}-e_{j} | {i}\neq{j}, \quad  i,j=1,...,n+1 \}$,
\end{center}
respectively.\\
\indent  Let $\Pi$ be a fundamental system in root system $\Phi$ of type $A_{n}$ and $\Phi^{+}$ be the 
corresponding positive system. Then $W(A_{n})$ is called the {\textit{ Weyl group of type
$A_{n}$}} generated by the reflections $w_{r}$ for all $r\in\Phi$.\newline

\indent  The descent algebra  $\sum(W)$  of a finite Coxeter group $W$ had been discovered by Solomon in 1976. After this seminal paper [2], there have been huge developments  in the area of  descent algebras of Coxeter groups over the last four decades. One can see for these developments, for example [3], [4] and [6].

\indent  Let $W_{J}$ be the subgroup of $W$, where $J$ is any subset of $\Pi$.
$W_{J}$ is called a standard parabolic subgroup of $W$.  Let $X_{J}$  be the
set of  representatives of the  cosets of $wW_{J}$ in $W$. Furthermore, $X_{J}^{-1}=\{w^{-1} | w\in{X_{J}}\}$ is a set of  representatives for the right cosets of  $W_{J}$ in $W$.
Then $X_{JK}=X_{J}^{-1}\bigcap{X_{K}}$ is a set of distinguished double coset representatives of $W_{J}wW_{K}$  in $W$.\newline

\indent The following remarkable theorem has been proved by
Solomon [2].\newline

\noindent\textbf{Theorem 1.1} {\textit{For every subset $J$ of   $\Pi$, let
$x_{J}=\sum_{d\in{X_{J}}}{d}$. Then 
\begin{center}
$x_Jx_K=\sum_{L\subseteq{K}}a_{JKL}x_L$,
\end{center}
where $a_{JKL}$  is the number of elements $w\in{X_{JK}}$ such that $w^{-1}(J)\cap{K}=L$.}}

\bigskip

\indent It follows that the set of all $x_{J}$ form a basis for a
noncommutative algebra $\sum(W)$ over the field of rationals. This algebra is
called {\textit{ the descent algebra (or Solomon algebra) of Coxeter groups $W$}} [3].\newline

\indent The descent algebra has been reconstructed by Bergeron et al. in 1992 [3]. In their paper, they have given a relatively elegant construction for the Solomon's descent algebra. Additionally, they have introduced the parabolic Burnside ring of the associated Coxeter group. Moreover, there is a natural homomorphism from the descent algebra $\sum(W)$ into the parabolic Burnside ring of the associated Coxeter group.\\ \\
\indent Clearly, since  the Weyl groups are also finite Coxeter groups then the construction of their descent algebra is valid for Weyl groups.\\
\indent  In other respects, the descent algebra associated with the  Weyl groups of type $A_{n}$ has been intensively studied. For example,  Garsia and Reutenauer [7] have given an  analysis of the descent algebra of the symmetric group. Besides, Atkinson [5] has determined the Loewy length of the descent algebra of the symmetric group.\\

\indent Now, for our purpose,  we  need  to define an equivalence relation on the set of all $x_J$, where $J\subseteq{\Pi}$. By this way, we can  obtain the equivalence classes in order to form the basis for a  new descent algebra.  Because of this, we now  give the following results.\\

\noindent\textbf{Definition 1.2} For $J,K\subseteq{\Pi}$, we write $J\sim{K}$ whenever $w(J)=K$ for
some $w\in{W}$, that is, $J$ and $K$ are {\textit{conjugate}} [3].\newline

This concept is used by Solomon [2] for groups as follows.\newline

\noindent\textbf{Lemma 1.3} {\textit{If $J,K\subseteq{\Pi}$ then $J\sim{K}$ if
and only if $W_{J}$ and $W_{K}$ are conjugate, that is, $W_{K}=wW_{J}w^{-1}$
for some $w\in{W}$.}}\newline

\indent Now, we define a relation on the set of all  $x_{J}$ in order to form a basis
for a new descent algebra of Weyl groups of type $A_{n}$  as follows:
\begin{eqnarray}
x_{J}\sim{x_{K}}\Leftrightarrow{J\sim{K}}\Leftrightarrow{K=w(J),\quad
for\quad w\in{W}}.
\end{eqnarray}

\noindent It is clear that  it is an equivalence relation. By means of this equivalence relation, for $J\subseteq{\Pi}$, the equivalence class
which contains $x_{J}$ is
\begin{align}
\lbrack x_{J}]  & =\{x_{K}|x_{K}\sim{x_{J}},K\subseteq\Pi\}\nonumber\\
& =\{x_{K}|K\sim{J}\}\nonumber\\
& =\{x_{w_{(J)}}|w(J)\subseteq{\Pi},w\in{W}\}.\nonumber
\end{align}
\\
\indent The main result of our paper is the following:\\

\noindent\textbf{Theorem A.} {\textit{Let $J\subseteq{\Pi}$, then
\begin{center}
$\sum_{W}(A_{n})=Sp\{ [x_J] | J\subseteq{\Pi}\}$
\end{center}
 is a new descent algebra of Weyl groups of type $A_{n}$.  }}\\

\section{Proof of Theorem A}
\bigskip
In this section, we shall obtain some significant results to prove Theorem A, and so we shall complete the construction. Furthermore, we shall obtain substantial results stating in Section 1.\\
\indent In that case, recall that the set $\{x_J | J\subseteq{\Pi}\} $ is a basis for the Solomon's descent algebra [2,3]. Then it is obvious that the set $ \{ [x_J] | J\subseteq{\Pi}\}$
is  linearly independence, and also a basis for $\sum_{W}(A_{n})$. Therefore $\sum_{W}(A_{n})$  is a $\mathbb{Q}%
$-space, which is  generated by all different equivalence classes $[x_{J}]$, where $J\subseteq{\Pi}$.\\ \\
\indent Besides, this $\mathbb{Q}%
$-space  $\sum_{W}(A_{n})$ must have a ring structure  for our purpose. Because of this, we now define a ring multiplication as follows:\\

\indent Let  $J,K\subseteq{\Pi}$, then
\begin{eqnarray}
[x_{J}][x_{K}]=[x_{J}x_{K}]=\sum_{L\subseteq{K}}a_{JKL}[x_{L}]
\end{eqnarray}
\newline where the $a_{JKL}$'s are defined as Section 1.\newline

\indent The following theorem shows that this multiplication is
well-defined.\newline

\noindent\textbf{Theorem 2.1} {\textit{If $J,J^{^{\prime}},K,K^{^{\prime}%
}\subseteq{\Pi}$, $J\sim{J^{^{\prime}}}$ and $K\sim{K^{^{\prime}}}$ then
\begin{align}
[x_{J}][x_{K}]=[x_{J^{^{\prime}}}][x_{K^{^{\prime}}}].\nonumber
\end{align}
\newline}}
\indent Now,  for the proof of this theorem,  we need to state  the following crucial  lemmas.\newline

\noindent\textbf{Lemma 2.2} {\textit{If $J,J^{^{\prime}},K,K^{^{\prime}%
}\subseteq{\Pi}$, $J\sim{J^{^{\prime}}}$ and $K\sim{K^{^{\prime}}}$ then
\begin{align}
|W_{J}\backslash{W}/W_{K}|=|W_{J^{^{\prime}}}\backslash{W}/W_{K^{^{\prime}}%
}|.\nonumber
\end{align}
\newline Proof.}} If $J\sim{J^{^{\prime}}}$ and $K\sim{K^{^{\prime}}}$ then  there
exist at least $w\in{W}$ and $\sigma\in{W}$ such that $J^{^{\prime}}=w(J)$ and
$K^{^{\prime}}=\sigma(K)$, respectively. Furthermore, by Lemma 1.3,  we can write
$W_{J^{^{\prime}}}=wW_{J}w^{-1}$ and $W_{K^{^{\prime}}}={\sigma}W_{J}{\sigma
}^{-1}$.\newline

\indent  Now we define a map as follows to achieve our goal:
\begin{center}
$\Psi:W_{J}\backslash{W}/W_{K}\rightarrow W_{J^{^{\prime}}}\backslash
{W}/W_{K^{^{\prime}}}$ ;   $\Psi(W_{J}zW_{K})=W_{J^{^{\prime}}}wz\sigma
^{-1}W_{K^{^{\prime}}}$
\end{center}

\noindent for $z\in{W}$.\newline \\
The map $\Psi$ is well-defined because, if\\
\begin{center}
$W_{J}z_{1}W_{K}=W_{J}z_{2}W_{K}$
\end{center}

\noindent then we obtain 
\begin{align}
\Psi(W_{J}z_{1}W_{K}) & =W_{J^{^{\prime}}}wz_{1}\sigma^{-1}W_{K^{^{\prime}}%
}\nonumber\\
& =wW_{J}z_{1}W_{K}\sigma^{-1}\nonumber\\
& =wW_{J}z_{2}W_{K}\sigma^{-1}\nonumber\\
& =\Psi(W_{J}z_{2}W_{K}).\nonumber
\end{align}
\noindent Obviously, if

\begin{center}
$\Psi(W_{J}z_{1}W_{K})=\Psi(W_{J}z_{2}W_{K})$
\end{center}

\noindent then we have 

\begin{center}
$W_{J}z_{1}W_{K}=W_{J}z_{2}W_{K}.$
\end{center}

\noindent Also, if we consider ${W_{J^{^{\prime}}}vW_{K^{^{\prime}}}}\in{ W_{J^{^{\prime}}}\backslash
{W}/W_{K^{^{\prime}}}}$ then there exist ${W_{J}w^{-1}v{\sigma}W_{K}}\in{W_{J}\backslash{W}/W_{K}}$ such that

\begin{center}
$\Psi(W_{J}w^{-1}v{\sigma}W_{K})=W_{J^{^{\prime}}}vW_{K^{^{\prime}}}.$
\end{center}

\indent So, this means that $\Psi$ is a bijective. Thus, the sets $W_{J}%
\backslash{W}/W_{K}$ and $W_{J^{^{\prime}}}\backslash{W}/W_{K^{^{\prime}}}$
have the same cardinality. This completes the proof. \quad  \quad  \quad \quad \quad \quad \quad  \quad  \quad  \quad \quad    \quad  $\Box$\newline\newline
\indent This lemma immediately implies the following corollary.
\newline

\noindent\textbf{Corollary 2.3} {\textit{For $J,J^{^{\prime}},K,K^{^{\prime}}\subseteq{\Pi}$, if
$J\sim{J^{^{\prime}}}$ and $K\sim{K^{^{\prime}}}$ then
\begin{align}
|X_{JK}|=|X_{J^{^{\prime}}K^{^{\prime}}}|.\nonumber
\end{align}
}}\\ 
\indent Also by Lemma 2.2 and Corollary 2.3,  we can see that there is
one-to-one correspondence between the elements of $X_{JK}$ and $X_{J^{^{\prime
}}K^{^{\prime}}}$. Because of this reason, we can obtain
\begin{center}
${W_{J^{^{\prime}}}fW_{K^{^{\prime}}}}= {W_{J^{^{\prime}}}wd{\sigma}%
^{-1}W_{K^{^{\prime}}}}=wW_{J}dW_{K}{\sigma}^{-1}.$
\end{center}

\noindent Accordingly, we can say that  there exist $w_{J}\in{W_{J}}$ and $w_{K}\in{W_{K}}$ such that
$f=ww_{J}dw_{K}{\sigma}^{-1}$.\newline\\
\indent The following theorem is very essential for our costruction.\\

\noindent\textbf{Theorem 2.4}{\textit{ Let ${J,K}\subseteq{\Pi}$ and $w\in{X_{JK}}.$ Then
\begin{center}
${W_J}\bigcap{wW_K{w^{-1}}}=W_L,$
\end{center}
where $L={J}\cap{w(K)}$ }}[3]. \\

\indent Now, we can give the following lemma.\\

\noindent\textbf{Lemma 2.5} {\textit{For $J,J^{^{\prime}},K,K^{^{\prime}%
}\subseteq{\Pi}$, $J\sim{J^{^{\prime}}}$ and $K\sim{K^{^{\prime}}}$, if
$d\in{X_{JK}}$ and $f\in{X_{{J^{^{\prime}}}{K^{^{\prime}}}}}$ then
\begin{align}
[x_{{f^{-1}}({J^{^{\prime}}})\cap{K^{^{\prime}}}}]=[x_{{d^{-1}}(J)\cap{K}%
}].\nonumber
\end{align}\\
Proof.}} By Lemma 1.3, if $J\sim{J^{^{\prime}}}$ and $K\sim{K^{^{\prime}}}$ then  there exist at
least $w\in{W}$ and $\sigma\in{W}$ such that $W_{J^{^{\prime}}}=wW_{J}w^{-1}$
and $W_{K^{^{\prime}}}={\sigma}W_{J}{\sigma}^{-1}$, respectively. And also  by Theorem 2.4, we know that the group  $W_{{f^{-1}}({J^{^{\prime}}})\cap{K^{^{\prime}}}}$ is equal to $f^{-1}W_{J^{^{\prime}}%
}f\bigcap{W_{K^{^{\prime}}}}$. Then
\begin{align}
W_{{f^{-1}}({J^{^{\prime}}})\cap{K^{^{\prime}}}} & =f^{-1}W_{J^{^{\prime}}%
}f\bigcap{W_{K^{^{\prime}}}}  \nonumber\\
& =\sigma{w_{K}}^{-1}d^{-1}{w_{J}}^{-1}w^{-1}wW_{J}w^{-1}ww_{J}dw_{K}\sigma
^{-1}\bigcap{{\sigma}W_{K}{\sigma^{-1}}}\nonumber\\
& =\sigma{w_{K}}^{-1}(d^{-1}W_{J}d\bigcap{W_{K}})w_{K}{\sigma^{-1}}\nonumber\\
& =\sigma{w_{K}}^{-1}W_{{d^{-1}}(J)\cap{K}}({\sigma}{w_{K}}^{-1}%
)^{-1}.\nonumber
\end{align}

\indent Note that, since ${\sigma}{w_{K}}^{-1}\in{W}$ then we obtain  $W_{{f^{-1}%
}({J^{^{\prime}}})\cap{K^{^{\prime}}}}\sim{W_{{d^{-1}}(J)\cap{K}}}$. This fact help us to
obtain  ${f^{-1}}({J^{^{\prime}}})\cap{K^{^{\prime}}}\sim{{d^{-1}}(J)\cap{K}}$, so
this means that 
\begin{center}
$[x_{{f^{-1}}({J^{^{\prime}}})\cap{K^{^{\prime}}}}]=[x_{{d^{-1}%
}(J)\cap{K}}]$.
\end{center}
This completes the proof. \quad \quad \quad \quad \quad \quad \quad \quad \quad \quad \quad \quad \quad \quad\quad \quad \quad\quad\quad \quad \quad \quad \quad $\Box$\newline

\indent We now have all the ingredients to give the proof of Theorem 2.1.\\ 

\noindent\textbf{Proof of  Theorem 2.1}.\newline

\noindent  Let $J,J^{^{\prime}},K,K^{^{\prime}%
}\subseteq{\Pi}$, $J\sim{J^{^{\prime}}}$ and $K\sim{K^{^{\prime}}}$. Then
\begin{align}
[x_{J^{^{\prime}}}][x_{K^{^{\prime}}}] & =[x_{J^{^{\prime}}}x_{K^{^{\prime}}%
}]\nonumber
\end{align}
\begin{align}
& =\sum_{f\in{X_{{J^{^{\prime}}}{K^{^{\prime}}}}}}[x_{{f^{-1}}({J^{^{\prime}}%
})\cap{K^{^{\prime}}}}]\nonumber\\
& =\sum_{d\in{X_{JK}}}[x_{{d^{-1}}(J)\cap{K}}]\nonumber\\
& =[x_{J}][x_{K}],\nonumber
\end{align}
and  Theorem 2.1 is proved. \quad \quad \quad \quad \quad \quad \quad \quad \quad \quad \quad \quad \quad \quad \quad  \quad \quad \quad \quad \quad  \quad \quad  $\Box$\\  \\
\indent Hence, the $\mathbb{Q}%
$-space  $\sum_{W}%
(A_{n})=Sp\{[x_{J}]|J\subseteq{\Pi}\}$ is a ring with identity $[x_{\Pi}]$. Obviously, this means that $\sum_{W}%
(A_{n})$ is an algebra over the field of rationals. 
Thus, the construction of the algebraic structure has been completed and also by means of this we complete the proof of Theorem A. 
To sum up, $\sum_{W}(A_{n})$  is a new descent algebra of Weyl groups of type $A_n$.\\

\indent Now, we want to determine two basic properties of this new descent algebra $\sum_{W}(A_{n})$, that is, this new descent algebra is commutative  and also semi-simple.\\

\indent By [2,3,4], it is known that  the radical of the descent algebra of Weyl groups is spanned by all differences $x_{J}-x_{K}$, where $J$ and $K$ are conjugate subsets of $\Pi$, that is to say,
\begin{center}
$Rad\sum(W)=Sp\{x_J-x_K | J\sim{K}, \quad J,K\subseteq{\Pi}\}.$
\end{center}
\bigskip
 Thus, by using this definition, we can obtain the radical of the new descent algebra of Weyl groups of type $A_{n}$ as follows:
\begin{center}
$Rad\sum_{W}(A_{n})=\{{0}\}.$
\end{center}
 
\bigskip
\indent Clearly, the next corollary is  the consequence of this result.\\

\noindent\textbf{Corollary 2.6} {\textit{The new descent algebra $\sum
_{W}(A_{n})$ of Weyl groups of type $A_{n}$ is semi-simple.}}\newline

\indent The other property of this algebra is given as follows.\newline

\noindent\textbf{Proposition 2.7} {\textit{The new descent algebra $\sum
_{W}(A_{n})$ of Weyl groups of type $A_{n}$ is commutative.}}\newline

\noindent{\textit{Proof.}} Let $J,K\subseteq{\Pi}$. Then
\begin{align}
[x_{J}][x_{K}] & =\sum_{d\in{X_{JK}}}[x_{{d^{-1}}(J)\cap{K}}]\nonumber
\end{align}
\begin{align}
& =\sum_{d^{-1}\in{X_{KJ}}}[x_{{d(K)\cap{J}}}]\nonumber\\
& =\sum_{f\in{X_{KJ}}}[x_{{f^{-1}}(K)\cap{J}}]\nonumber\\
& =[x_{K}][x_{J}].\nonumber
\end{align}
So, Proposition 2.7 is proved.
\quad \quad \quad \quad \quad \quad \quad \quad \quad \quad \quad \quad \quad \quad \quad \quad \quad \quad \quad \quad \quad $\Box$
\newline
\
\bigskip
\\ 
\noindent\textbf{Remark 2.8}{ We have seen that the new descent algebra $\sum_{W}(A_{n})$ of Weyl groups of type $A_{n}$ generated by equivalence classes is different from the Solomon algebra  $\sum(W)$ of Weyl groups of type $A_{n}$, since Solomon  algebra is non-commutative and not semi-simple. } \\

\indent The following example illustrates Theorem A.\\

\noindent\textbf{Example 2.9}  Let us give the multiplication table for the new descent algebra $\sum_{W}(A_{3})$, where the fundamental system of type $A_{3}$ is 
\begin{center}
 $\Pi =\{e_{1}-e_{2}, e_{2}-e_{3}, e_{3}-e_{4}\}$.
\end{center}
\bigskip

\indent By Theorem A, we know that the new descent algebra of Weyl group type  $A_{3}$ is generated by the equivalence classes arising from the  equivalence relation (1) defined in Section 1. Thus, by using this equivalence relation and also with keeping in mind the distinguished basis of the Solomon's descent algebra of Weyl group of type $A_{3}$, which is  given by the sums of the distinguished coset representatives of parabolic subgroups of $W(A_3)$,  then  we can obtain the generators of this new descent algebra as follows:\\

 \indent $[x_J], \quad  where \quad   {J={\Pi}}$,\\
 \indent $[x_K]=[x_L], \quad  where \quad {K=\{e_{1}-e_{2}, e_{2}-e_{3}\}} \quad or \\
 \indent \quad \quad \quad \quad \quad \quad \quad \quad \quad \quad  {L=\{e_{2}-e_{3}, e_{3}-e_{4}\}}$,\\
\indent $[x_M],\quad where  \quad {M=\{e_{1}-e_{2}, e_{3}-e_{4}\}}$,\\
 \indent $[x_N]=[x_P]=[x_Q],  \quad where \quad {N=\{e_{1}-e_{2}\}}     \quad or\\
\indent \quad \quad \quad \quad \quad \quad \quad \quad \quad \quad \quad \quad  \quad  {P=\{e_{2}-e_{3}\}} \quad or \\
\indent \quad \quad \quad \quad \quad \quad \quad \quad  \quad \quad \quad \quad  \quad {Q=\{e_{3}-e_{4}\}}$,\\
 \indent $[x_R], \quad where \quad {R={\emptyset}}$.\\

If we denote these elements by $x_{{\lambda}_1}$,  $x_{{\lambda}_2}$,   $x_{{\lambda}_3}$,  $x_{{\lambda}_4}$ and   $x_{{\lambda}_5}$,  respectively,  where ${{\lambda}_1}=J$,  ${{\lambda}_2}=K${\textit{ or}} L, ${{\lambda}_3}=M$,  ${{\lambda}_4}=N$ {\textit{or}} P {\textit{or}}  Q,  ${{\lambda}_5}=R$, then we get Table 1 as follows by using ring multiplication operation (2).\\  

\begin{table}
\begin{center}

\begin{tabular}{ |  c  |  c  |  c  |  c  |  c  |   c   |  }
\cline{2-6}  
\multicolumn{1}{c|}{} &$x_{{\lambda}_1}$&$x_{{\lambda}_2}$& $x_{{\lambda}_3}$& $x_{{\lambda}_4}$& $x_{{\lambda}_5}$\\
\hline
$x_{{\lambda}_1}$ & $x_{{\lambda}_1}$ & $x_{{\lambda}_2}$ & $x_{{\lambda}_3}$&$x_{{\lambda}_4}$&$x_{{\lambda}_5}$ \\
\hline
$x_{{\lambda}_2}$ &$x_{{\lambda}_2}$ & $x_{{\lambda}_2}+x_{{\lambda}_4}$  &$2{x_{{\lambda}_4}}$ &$2{x_{{\lambda}_4}}+x_{{\lambda}_5}$&$4{x_{{\lambda}_5}}$\\
\hline
$x_{{\lambda}_3}$ & $x_{{\lambda}_3}$ &$2{x_{{\lambda}_4}}$ & $2{x_{{\lambda}_3}}+x_{{\lambda}_5}$ &$2{x_{{\lambda}_4}}+2{x_{{\lambda}_5}}$&$6{x_{{\lambda}_5}}$ \\
\hline
$x_{{\lambda}_4}$ & $x_{{\lambda}_4}$ &$2{x_{{\lambda}_4}}+x_{{\lambda}_5}$ & $2{x_{{\lambda}_4}}+2{x_{{\lambda}_5}}$ &$2{x_{{\lambda}_4}}+5{x_{{\lambda}_5}}$&$12{x_{{\lambda}_5}}$\\
\hline
$x_{{\lambda}_5}$ & $x_{{\lambda}_5}$ & $4{x_{{\lambda}_5}}$ & $6{x_{{\lambda}_5}}$&$12{x_{{\lambda}_5}}$&$24{x_{{\lambda}_5}}$\\
\hline
\end{tabular}
\end{center}
\caption{ The ring multiplication table of $\sum_{W}(A_{3})$.}
\label{Table1}
\end{table}

 \indent This example  shows that the new commutative descent  algebra of Weyl group of type    $A_{3}$  has less elements than the  Solomon's  algebra associated with $W(A_3)$. \\

\indent  In other respects, the descent algebra is closely related to the subring of the Burnside ring $B(W)$. This ring is called the {\textit{ parabolic Burnside ring}}  and denoted by $PB(W)$. Furthermore, the  parabolic Burnside ring is spanned by the permutation representations $W/W_{J}$, where the $W_J$ are the parabolic subgroups of $W$. And also, the multiplication for $ W/W_{J}$ and   $ W/W_{K}$ is defined by
\begin{align}
 {W/W_{J}}\times{W/W_{K}}=\sum_{L\subseteq{K}}a_{JKL}{W/W_{L}}
\end{align}
where the $a_{JKL}$'s are defined as Section 1 [3].\\
\indent Now we can state the following proposition.\\

\noindent\textbf{Proposition 2.10} {\textit{ Let $\phi:\sum_{W}(A_{n})\rightarrow PB(W(A_n))$ be the map given by $\phi([x_J])=W/W_{J}$. This is an isomorphism.}}\\

\noindent{\textit{Proof.}} This mapping is well-defined. In fact, if we take $[x_J]=[x_K]$ for $J,K \subseteq{\Pi}$  then we can obtain $W/W_{J}\cong{W/W_{K}}$ because of $J\sim{K}$ (see [3]). So, this implies that $\phi([x_J])=\phi([x_K])$. Additionally, one can easily verify that $\phi$ is a bijective. To show that $\phi$ is an isomorphism, it is enough to show that $\phi([x_J][x_K])={\phi([x_J])}{\phi([x_K])}$.  By using the equations (2) and (3), we obtain
\begin{align}
\phi([x_J][x_K]) & =\phi([{x_J}{x_K}])\nonumber\\
& =\phi(\sum_{L\subseteq{K}}a_{JKL}[x_{L}])\nonumber\\
& =\sum_{L\subseteq{K}}a_{JKL}\phi([x_{L}])\nonumber\\
& =\sum_{L\subseteq{K}}a_{JKL}{W/W_{L}}\nonumber\\
&= {W/W_{J}}\times{W/W_{K}}\nonumber\\
&={\phi([x_J])}{\phi([x_K])}\nonumber,
\end{align}
\noindent and thus, Proposition 2.10 is proved. \quad \quad \quad  \quad \quad\quad   \quad \quad \quad  \quad \quad \quad \quad  \quad \quad \quad \quad  \quad $\Box$\\

\indent Let ${\chi}_J$ be the permutation character of $ W(A_n)$ acting on the left cosets of $W_J$  and let $G(W)$  be the $\mathbb{Z}$-module generated by all ${\chi}_J$ [3,4]. The following proposition was proved in [3].\\

\noindent\textbf{Proposition 2.11} {\textit{ The assignment $W/W_{J}\mapsto{\chi_J}$ defines an isomorphism $\Theta$ from $PB(W)$ to the $G(W)$.}}\\

\indent From this point of view,  there is an isomorphism between the new descent algebra $\sum_{W}(A_{n})$  and  $G(W(A_n))$. This isomorphism is stated as follows:\\

\noindent\textbf{Proposition 2.12} {\textit{ The linear map $\psi$  defined by the images $\psi([x_J])={\chi}_J$ is an isomorphism from  the new descent algebra $\sum_{W}(A_{n})$ to the  $G(W(A_n))$.}}\\

\noindent{\textit{Proof.}}   The linear map $\psi$ is a composition of $\Theta$ and  $\phi$  defined in Proposition 2.11 and 2.10, respectively. This completes the proof.  \quad \quad \quad \quad \quad \quad \quad  \quad \quad \quad \quad \quad \quad$\Box$\\

\indent Thus, we may identify the new descent algebra with this ring of class functions.

\end{document}